\documentclass{article}
\usepackage{amsmath}

\input{tcilatex}
\begin{document}

\begin{center}
\textbf{New Estimations for Sturm\textbf{-}Liouville Problems in Difference
Equations}

\bigskip

Erdal BAS, \ Ramazan OZARSLAN

\bigskip

\bigskip \textit{Department of Mathematics, Faculty of Science, Firat
University, Elazig, 23119, Turkey}

\textit{e-mail: erdalmat@yahoo.com, }ozarslanramazan@gmail.com\bigskip
\end{center}

\noindent \textbf{Abstract}

{\footnotesize In this paper, Sturm-Liouville problem for difference
equations is considered with potential function }$q\left( n\right) .$%
{\footnotesize \ The representations of solutions are obtained by variation
of parameters method. These solutions are proved, using summation by parts.
Also, estimation\ of asymptotic expansion of the solutions are established.}

\bigskip 

\noindent \textbf{AMS:} 39A13, 34B18.\bigskip 

\bigskip

\noindent \textbf{1. Introduction}

Sturm-Liouville operators have been studied for a number of years. Firstly,
Sturm Liouville problem developed in a number of articles published by these
authors in 1836 and 1837. It is known that the spectral characteristics are
spectra, spectral functions, scattering data, norming constants, etc. The
representation of solution of Sturm-Liouville problem and asymptotic
formulas for eigenfunctions have been obtained by Levitan and Sargsjan $%
\left[ 6\right] $. The following problem%
\begin{equation*}
Ly\left( t\right) =-\frac{d^{2}y}{dt^{2}}+q\left( t\right) y%
\begin{array}{c}
=%
\end{array}%
\lambda y,
\end{equation*}%
\begin{eqnarray*}
&&y\left( 0\right) \cos \alpha +y^{\prime }\left( 0\right) \sin \alpha 
\begin{array}{c}
=%
\end{array}%
0 \\
&&y\left( 1\right) \cos \beta +y^{\prime }\left( 1\right) \sin \beta 
\begin{array}{c}
=%
\end{array}%
0
\end{eqnarray*}%
is called Sturm-Liouville problem in differential equation. Differential
equations are related to difference equations closely.

In general, it is known that difference equations related to recursive
relations for a long time. Its actual development appeared by being able to
compared to differential equations. Difference equations have many
application areas which are problems of physics, mathematics and
engineering, vibrating string, economy, population, actueria and logistics,
etc.

As a natural result of comparing the difference equations to differential
equations, the theory of linear difference equations have begun to appear in
a similar way to the theory of differential equations. Basic theory of
linear difference equations improved by De Moivre, Euler, Lagrange, Laplace
et al $\left[ 13\right] $. Hereafter, Hartman $\left[ 17\right] $,
Ahlbrandt-Hooker $\left[ 16\right] $, Peterson-Kelley $\left[ 1\right] $,
Agarwal $\left[ 2\right] $, Jirari $\left[ 5\right] $, Bender-Orszag $\left[
7\right] $, Goldberg $\left[ 12\right] $, Elaydi $\left[ 15\right] $,
Lakshmikantham, Trigiante and Peterson $\left[ 13\right] $, Mickens $\left[
14\right] $ and they have contributed to linear difference equations with
publications and books.

Especially, in recent years, Sturm-Liouville difference equation has seen a
great interest and a lot of study has published, but there are still a lot
of things to develop about this subject. Peterson, Kelley, Agarwal, Bender
gave a place this subject in studies, also, Jirari studied pecularly
Sturm-Liouville difference equation in his studies.

Atkinson studied discrete and continuous boundary value problems in his book
and also he considered self-adjoint second-order difference equations.
Jirari contributed Atkinson's study by Second Order Sturm-Liouville
Difference Equations and Orthogonal Polynomials in his thesis.

Hinton and Lewis $\left[ 9\right] $, Clark $\left[ 19\right] $, Shi and Wu $%
\left[ 8\right] ,$ Shi and Sun $\left[ 10\right] ,$ Shi and Chen $\left[ 18%
\right] \ $and Hilscher $\left[ 20\right] $ investigated spectral analysis
of second order difference equations and operators. Wang and Shi $\left[ 11%
\right] ,$ Ji and Yang $\left[ 21\right] ,$ Sun and Shi $\left[ 22\right] $
studied eigenvalues of second order difference equations.

In the most of studies for Sturm-Liouville difference equations, $q\left(
n\right) $ is considered as a real number, solutions and spectral properties
for Sturm-Liouville difference equation are investigated by Hamiltoni
systems. The following $(1.1)-(1.3)$ problem $\left[ 1\right] ,$ $\left[ 5%
\right] $%
\begin{eqnarray}
\Delta \left( \Delta p\left( n-1\right) x\left( n-1\right) \right) +q\left(
n\right) x\left( n\right) +\lambda r\left( n\right) x\left( n\right) &=&0,%
\text{ }n=a,...,b  \TCItag{1.1} \\
x\left( a-1\right) +hx\left( a\right) &=&0,  \TCItag{1.2} \\
x\left( b+1\right) +kx\left( b\right) &=&0.  \TCItag{1.3}
\end{eqnarray}%
is called Sturm-Liouville problem in difference equations.

Our study is organized as follows. Fundamental theorems and definitions in
Section 2, representations of solutions with two different initial
conditions in Section 3 and asymptotic behavior of eigenfunctions are given
in Section 4.

\bigskip

\noindent \textbf{2. Preliminaries}

\bigskip

\noindent \textbf{Definition.2.1. }$\left[ 1\right] $ The matrix of
Casoratian is given by%
\begin{equation*}
w\left( n\right) =\left( 
\begin{array}{cccc}
x_{1}\left( n\right) & x_{2}\left( n\right) & ... & x_{r}\left( n\right) \\ 
x_{1}\left( n+1\right) & x_{2}\left( n+1\right) & ... & x_{r}\left(
n+1\right) \\ 
\vdots & \vdots &  & \vdots \\ 
x_{1}\left( n+r-1\right) & x_{2}\left( n+r-1\right) & ... & x_{r}\left(
n+r-1\right)%
\end{array}%
\right)
\end{equation*}%
where $x_{1}\left( n\right) ,x_{2}\left( n\right) ,...,x_{r}\left( n\right) $
are given functions. The determinant%
\begin{equation*}
W\left( n\right) =\det w\left( n\right)
\end{equation*}%
is called Casoratian. Note that, it is similar to Wronskian determinant.

\bigskip

\noindent \textbf{Theorem.2.1. }$\left[ 5\right] $ (Wronskian-Type Identity)
Let $y$ and $z$ be a solutions of $\left( 1.1\right) .$ Then, for $a\leq
n\leq b$%
\begin{eqnarray}
W\left[ y,z\right] \left( n\right) &=&p\left( n-1\right) \left[ y\left(
n\right) \Delta z\left( n-1\right) -z\left( n\right) \Delta y\left(
n-1\right) \right]  \TCItag{2.1} \\
&=&-p\left( n-1\right) \left[ y\left( n\right) z\left( n-1\right) -y\left(
n-1\right) z\left( n\right) \right]  \notag
\end{eqnarray}%
is a constant (In particular equal to $W\left[ y,z\right] \left( a\right) $).

\bigskip

\noindent \textbf{Definition.2.2. }$\left[ 5\right] $ Let's express $\left(
1.1\right) $ Sturm-Liouville equation as follows,%
\begin{equation}
Lx\left( n\right) =-\lambda x\left( n\right) ,\text{ }n\in \left[ a,b\right]
,  \tag{2.2}
\end{equation}%
with initial conditions%
\begin{equation}
\cos \alpha x\left( a\right) -\sin \alpha \left( p\left( a\right) \nabla
x\left( a\right) \right) =0,  \tag{2.3}
\end{equation}%
where $0<\alpha ,$ $\beta <\pi ,$ $\nabla $ is the backward difference
operator, $\nabla x\left( n\right) =x\left( n\right) -x\left( n-1\right) ,$
are equivalent to%
\begin{equation}
x\left( a-1\right) +\left( \frac{\cot \alpha }{p\left( a\right) }-1\right)
x\left( a\right) =0,  \tag{2.4}
\end{equation}%
in other words,%
\begin{equation}
x\left( a-1\right) +hx\left( a\right) =0,  \tag{2.5}
\end{equation}%
where $L$ is self-adjoint Sturm-Liouville operator and $h$ is real number. $%
\left( 2.2\right) -\left( 2.5\right) $ initial value problem is called
Sturm-Liouville problem.

\bigskip

\noindent \textbf{Theorem.2.2. }$\left[ 1\right] $ $\left( \text{\textbf{%
Summation by parts}}\right) $ If $m<n$, then%
\begin{equation}
\sum_{k=m}^{n-1}x\left( k\right) \Delta y\left( k\right) =\left[ x\left(
k\right) y\left( k\right) \right] _{m}^{n}-\sum_{k=m}^{n-1}\Delta x\left(
k\right) y\left( k+1\right) .  \tag{2.6}
\end{equation}

\noindent \textbf{Theorem.2.4}. $\left[ 1\right] $ \ If $y_{n}$ is an
indefinite sum of $x_{n}$, then%
\begin{equation}
\sum_{k=m}^{n-1}y\left( k\right) =x\left( n\right) -x\left( m\right) . 
\tag{2.7}
\end{equation}

\bigskip

\noindent \textbf{Theorem.2.5.} $\left[ 1\right] ,\left[ 3\right] $ $\left( 
\text{\textbf{Annihilator method}}\right) $ Suppose that $x\left( n\right) $
solves following difference equation, $E$ is the shift operator, $Ex\left(
n\right) =x\left( n+1\right) ,$%
\begin{equation*}
\left( E^{t}+p\left( t-1\right) E^{t-1}+...+p\left( 0\right) \right) x\left(
n\right) =r\left( n\right) ,
\end{equation*}%
and that $r\left( n\right) $ satisfies%
\begin{equation*}
\left( E^{m}+q\left( m-1\right) E^{m-1}+...+q\left( 0\right) \right) r\left(
n\right) =0.
\end{equation*}%
Then $x\left( n\right) $ satisfies%
\begin{equation*}
\left( E^{m}+q\left( m-1\right) E^{m-1}+...+q\left( 0\right) \right) \left(
E^{t}+p\left( t-1\right) E^{t-1}+...+p\left( 0\right) \right) x\left(
n\right) =0.
\end{equation*}

\bigskip

\textbf{3. Main Results}

In this paper, we are interested in the representations of solutions of the
Sturm-Liouville problem in difference equations with potential function $%
q\left( n\right) $ as follows,

\begin{equation}
\Delta ^{2}x\left( n-1\right) +q\left( n\right) x\left( n\right) +\lambda
x\left( n\right) =0,\text{ }n=a,...,b  \tag{3.1}
\end{equation}%
with initial conditions,%
\begin{equation}
x\left( a-1\right) +hx\left( a\right) =0,  \tag{3.2}
\end{equation}%
where $a,$ $b$ are finite integers with $a\geq 0,$ $a\leq b,$ $h$ is a real
number, $\Delta $ is the forward difference operator, $\Delta x\left(
n\right) =x\left( n+1\right) -x\left( n\right) ,$ $\lambda $ is the spectral
parameter, $q\left( n\right) $ is a real valued potential function for $n\in %
\left[ a,b\right] .$ In the general literature, potential function $q\left(
n\right) $ is taken as real number but we take it as a variable coefficients
in a similar way in Levitan and Sargsjan's study $\left[ 6\right] $. Levitan
and Sargsjan obtained the representation of solution of Sturm-Liouville
problem in differential equations and asymptotic formulas for
eigenfunctions. In analogous manner, we tried to obtain the representation
of solution of Sturm-Liouville problem in difference equations.

\noindent A self-adjoint difference operator corresponds to equation $\left(
3.1\right) $ noted by,%
\begin{equation*}
Lx\left( n\right) =\Delta ^{2}x\left( n-1\right) +q\left( n\right) x\left(
n\right) =-\lambda x\left( n\right) .
\end{equation*}

\noindent In $\ell ^{2}\left( a,b\right) ,$ the Hilbert space of sequences
of complex numbers $x\left( a\right) ,...,x\left( b\right) $ with the inner
product,%
\begin{equation*}
<x\left( n\right) ,y\left( n\right) >=\sum_{n=a}^{b}x\left( n\right) y\left(
n\right) ,
\end{equation*}

\noindent for every $x\in D_{L,}$ let define as follows%
\begin{equation*}
D_{L}=\left\{ x\left( n\right) \in \ell ^{2}\left( a,b\right) :Lx\left(
n\right) \in \ell ^{2}\left( a,b\right) ,\text{ }x\left( 0\right) =-h,\text{ 
}x\left( 1\right) =1\right\} .
\end{equation*}%
Hence, equation $\left( 3.1\right) $ can be written as follows%
\begin{equation*}
Lx\left( n\right) =-\lambda x\left( n\right) .
\end{equation*}%
At this section, we present the representation of solution of $\left(
3.1\right) -\left( 3.2\right) $ Sturm-Liouville problem by variation of
parameters method.

\bigskip

\noindent \textbf{Theorem 3.1. }Let define Sturm-Liouville problem in
difference equations as follows;%
\begin{equation}
Lx\left( n\right) =-\lambda x\left( n\right) ,  \tag{3.3}
\end{equation}%
\begin{equation}
x\left( 0\right) =-h,\text{ }x\left( 1\right) =1,  \tag{3.4}
\end{equation}%
then $\left( 3.3\right) -\left( 3.4\right) $ Sturm-Liouville problem has a
unique solution for $x\left( n\right) $ as%
\begin{eqnarray}
x\left( n,\lambda \right) &=&\left( \frac{2-h\left( 2q\left( 0\right)
-2+\lambda +\sqrt{\lambda \left( \lambda -4\right) }\right) }{2\sqrt{\lambda
\left( \lambda -4\right) }}\right) \left( \frac{2-\lambda +\sqrt{\lambda
\left( \lambda -4\right) }}{2}\right) ^{n}  \TCItag{3.5} \\
&&+\left( \frac{-2+h\left( 2q\left( 0\right) -2+\lambda -\sqrt{\lambda
\left( \lambda -4\right) }\right) }{2\sqrt{\lambda \left( \lambda -4\right) }%
}\right) \left( \frac{2-\lambda -\sqrt{\lambda \left( \lambda -4\right) }}{2}%
\right) ^{n}  \notag \\
&&-\sum_{i=0}^{n}\frac{q\left( i\right) x\left( i\right) \left( \frac{%
2-\lambda -\sqrt{\lambda \left( \lambda -4\right) }}{2}\right) ^{i}}{\sqrt{%
\lambda \left( \lambda -4\right) }}\left( \frac{2-\lambda +\sqrt{\lambda
\left( \lambda -4\right) }}{2}\right) ^{n}  \notag \\
&&+\sum_{i=0}^{n}\frac{q\left( i\right) x\left( i\right) \left( \frac{%
2-\lambda +\sqrt{\lambda \left( \lambda -4\right) }}{2}\right) ^{i}}{\sqrt{%
\lambda \left( \lambda -4\right) }}\left( \frac{2-\lambda -\sqrt{\lambda
\left( \lambda -4\right) }}{2}\right) ^{n}.  \notag
\end{eqnarray}%
Where $\sum\limits_{i=0}^{-1}.=0.$

\bigskip

\noindent \textbf{Proof.} If $x_{1}\left( n\right) $ and $x_{2}\left(
n\right) $ are linearly independent solution for homogen part of $\left(
3.1\right) ,$ then it is easily found by characteristic polynomial $\left[ 1%
\right] $%
\begin{equation}
x_{h}\left( n\right) =c_{1}x_{1}\left( n\right) +c_{2}x_{2}\left( n\right) ,
\tag{3.6}
\end{equation}%
\begin{equation}
x_{h}\left( n\right) =c_{1}\left( \frac{2-\lambda +\sqrt{\lambda \left(
\lambda -4\right) }}{2}\right) ^{n}+c_{2}\left( \frac{2-\lambda -\sqrt{%
\lambda \left( \lambda -4\right) }}{2}\right) ^{n}.  \tag{3.7}
\end{equation}%
Assume that $\mid \lambda -2\mid <2$ for existence of eigenvalues of $\left(
3.7\right) .$ By variation of parameters method $\left[ 1\right] ,$ $\left[ 7%
\right] $, we take 
\begin{equation}
x_{p}\left( n\right) =c_{1}\left( n\right) x_{1}\left( n\right) +c_{2}\left(
n\right) x_{2}\left( n\right) ,  \tag{3.8}
\end{equation}

\noindent If we perform necessary processes, we find $c_{1}$ and $c_{2}$ as
follows%
\begin{eqnarray}
c_{1}\left( n\right)  &=&\sum_{i=0}^{n}\frac{q\left( i\right) x\left(
i\right) x_{2}\left( i\right) }{W\left( x_{1}\left( i\right) ,x_{2}\left(
i\right) \right) },  \TCItag{3.9} \\
c_{2}\left( n\right)  &=&\sum_{i=0}^{n}\frac{-q\left( i\right) x\left(
i\right) x_{1}\left( i\right) }{W\left( x_{1}\left( i\right) ,x_{2}\left(
i\right) \right) }.  \notag
\end{eqnarray}%
Finally, we obtain general solution%
\begin{eqnarray*}
x\left( n,\lambda \right)  &=&c_{1}\left( \frac{2-\lambda +\sqrt{\lambda
\left( \lambda -4\right) }}{2}\right) ^{n}+c_{2}\left( \frac{2-\lambda -%
\sqrt{\lambda \left( \lambda -4\right) }}{2}\right) ^{n} \\
&&+\sum_{i=0}^{n}\frac{q\left( i\right) x\left( i\right) \left( \frac{%
2-\lambda -\sqrt{\lambda \left( \lambda -4\right) }}{2}\right) ^{i}}{W\left(
x_{1}\left( i\right) ,x_{2}\left( i\right) \right) }\left( \frac{2-\lambda +%
\sqrt{\lambda \left( \lambda -4\right) }}{2}\right) ^{n} \\
&&+\sum_{i=0}^{n}\frac{-q\left( i\right) x\left( i\right) \left( \frac{%
2-\lambda +\sqrt{\lambda \left( \lambda -4\right) }}{2}\right) ^{i}}{W\left(
x_{1}\left( i\right) ,x_{2}\left( i\right) \right) }\left( \frac{2-\lambda -%
\sqrt{\lambda \left( \lambda -4\right) }}{2}\right) ^{n}.
\end{eqnarray*}%
Where, $W$ Casoratian determinant is a constant by Theorem 2.1,%
\begin{equation*}
W\left( x_{1}\left( i\right) ,x_{2}\left( i\right) \right) =-\sqrt{\lambda
\left( \lambda -4\right) }.
\end{equation*}%
If we use the initial conditions $\left( 3.4\right) $, then we obtain the
representation of the solution of the Sturm-Liouville problem in difference
equations as follows,%
\begin{eqnarray}
x\left( n,\lambda \right)  &=&\left( \frac{2-h\left( 2q\left( 0\right)
-2+\lambda +\sqrt{\lambda \left( \lambda -4\right) }\right) }{2\sqrt{\lambda
\left( \lambda -4\right) }}\right) \left( \frac{2-\lambda +\sqrt{\lambda
\left( \lambda -4\right) }}{2}\right) ^{n}  \TCItag{3.10} \\
&&+\left( \frac{-2+h\left( 2q\left( 0\right) -2+\lambda -\sqrt{\lambda
\left( \lambda -4\right) }\right) }{2\sqrt{\lambda \left( \lambda -4\right) }%
}\right) \left( \frac{2-\lambda -\sqrt{\lambda \left( \lambda -4\right) }}{2}%
\right) ^{n}  \notag \\
&&-\sum_{i=0}^{n}\frac{q\left( i\right) x\left( i\right) \left( \frac{%
2-\lambda -\sqrt{\lambda \left( \lambda -4\right) }}{2}\right) ^{i}}{\sqrt{%
\lambda \left( \lambda -4\right) }}\left( \frac{2-\lambda +\sqrt{\lambda
\left( \lambda -4\right) }}{2}\right) ^{n}  \notag \\
&&+\sum_{i=0}^{n}\frac{q\left( i\right) x\left( i\right) \left( \frac{%
2-\lambda +\sqrt{\lambda \left( \lambda -4\right) }}{2}\right) ^{i}}{\sqrt{%
\lambda \left( \lambda -4\right) }}\left( \frac{2-\lambda -\sqrt{\lambda
\left( \lambda -4\right) }}{2}\right) ^{n}.  \notag
\end{eqnarray}

\bigskip

\noindent Now, let's show that $\left( 3.15\right) $ holds Sturm-Liouville
problem $\left( 3.3\right) -\left( 3.4\right) $. From $\left( 3.3\right) $%
\begin{equation}
q\left( n\right) x\left( n\right) =-\Delta ^{2}x\left( n-1\right) -\lambda
x\left( n\right) .  \tag{3.11}
\end{equation}%
First, let's take last two term in $\left( 3.10\right) $ and write equality $%
\left( 3.11\right) $ in place of $q\left( i\right) x\left( i\right) $.
Hence, we obtain%
\begin{equation}
-\sum_{i=0}^{n}\tfrac{\left[ -\Delta ^{2}x\left( i-1\right) -\lambda x\left(
i\right) \right] x_{2}\left( i\right) }{\sqrt{\lambda \left( \lambda
-4\right) }}x_{1}\left( n\right) =\sum_{i=0}^{n}\tfrac{\Delta ^{2}x\left(
i-1\right) x_{2}\left( i\right) }{\sqrt{\lambda \left( \lambda -4\right) }}%
x_{1}\left( n\right) +\sum_{i=0}^{n}\tfrac{\lambda x\left( i\right)
x_{2}\left( i\right) }{\sqrt{\lambda \left( \lambda -4\right) }}x_{1}\left(
n\right) ,  \tag{3.12}
\end{equation}%
\begin{equation}
\sum_{i=0}^{n}\tfrac{\left[ -\Delta ^{2}x\left( i-1\right) -\lambda x\left(
i\right) \right] x_{1}\left( i\right) }{\sqrt{\lambda \left( \lambda
-4\right) }}x_{2}\left( n\right) =-\sum_{i=0}^{n}\tfrac{\Delta ^{2}x\left(
i-1\right) x_{1}\left( i\right) }{\sqrt{\lambda \left( \lambda -4\right) }}%
x_{2}\left( n\right) -\sum_{i=0}^{n}\tfrac{\lambda x\left( i\right)
x_{1}\left( i\right) }{\sqrt{\lambda \left( \lambda -4\right) }}x_{2}\left(
n\right)   \tag{3.13}
\end{equation}%
Second, let's apply twice summation by parts method to first term at the
right hand side of equation $\left( 3.12\right) $ and $\left( 3.13\right) $
by Theorem 2.2, we obtain%
\begin{eqnarray}
-\sum_{i=0}^{n}\frac{\left[ -\Delta ^{2}x\left( i-1\right) -\lambda x\left(
i\right) \right] x_{2}\left( i\right) }{\sqrt{\lambda \left( \lambda
-4\right) }}x_{1}\left( n\right)  &=&[x_{2}\left( n+1\right) \Delta x\left(
n\right) -x_{2}\left( 0\right) \Delta x\left( -1\right)   \TCItag{3.14} \\
&&-\Delta x_{2}\left( n+1\right) x\left( n+1\right) +\Delta x_{2}\left(
0\right) x\left( 0\right)   \notag \\
&&-\Delta ^{2}x_{2}\left( -1\right) x\left( 0\right) ]\frac{x_{1}\left(
n\right) }{\sqrt{\lambda \left( \lambda -4\right) }}  \notag \\
&&+\sum\limits_{i=0}^{n}\frac{\left[ \Delta ^{2}x_{2}\left( i-1\right)
+\lambda x_{2}\left( i\right) \right] x\left( i\right) }{\sqrt{\lambda
\left( \lambda -4\right) }}x_{1}\left( n\right)   \notag
\end{eqnarray}%
and%
\begin{eqnarray}
\sum_{i=0}^{n}\frac{\left[ -\Delta ^{2}x\left( i-1\right) -\lambda x\left(
i\right) \right] x_{1}\left( i\right) }{\sqrt{\lambda \left( \lambda
-4\right) }}x_{2}\left( n\right)  &=&-[x_{1}\left( n+1\right) \Delta x\left(
n\right) -x_{1}\left( 0\right) \Delta x\left( -1\right)   \TCItag{3.15} \\
&&-\Delta x_{1}\left( n+1\right) x\left( n+1\right) +\Delta x_{1}\left(
0\right) x\left( 0\right)   \notag \\
&&-\Delta ^{2}x_{1}\left( -1\right) x\left( 0\right) ]\frac{x_{2}\left(
n\right) }{\sqrt{\lambda \left( \lambda -4\right) }}  \notag \\
&&-\sum\limits_{i=0}^{n}\frac{\left[ \Delta ^{2}x_{1}\left( i-1\right)
+\lambda x_{1}\left( i\right) \right] x\left( i\right) }{\sqrt{\lambda
\left( \lambda -4\right) }}x_{2}\left( n\right) .  \notag
\end{eqnarray}

\noindent Since $x_{1}$ and $x_{2}$ satisfy the homogen part of $\left(
3.3\right) $, sum expressions at the right hand side of $\left( 3.14\right) $
and $\left( 3.15\right) $ equal to zero and hence, $\left( 3.14\right) $ and 
$\left( 3.15\right) $ as follows respectively,%
\begin{eqnarray}
&&[x_{2}\left( n+1\right) \Delta x\left( n\right) -x_{2}\left( 0\right)
\Delta x\left( -1\right)   \TCItag{3.16} \\
&&-\Delta x_{2}\left( n+1\right) x\left( n+1\right) +\Delta x_{2}\left(
0\right) x\left( 0\right)   \notag \\
&&-\Delta ^{2}x_{2}\left( -1\right) x\left( 0\right) ]\frac{x_{1}\left(
n\right) }{\sqrt{\lambda \left( \lambda -4\right) }},  \notag
\end{eqnarray}

\begin{eqnarray}
&&-[x_{1}\left( n+1\right) \Delta x\left( n\right) -x_{1}\left( 0\right)
\Delta x\left( -1\right)   \TCItag{3.17} \\
&&-\Delta x_{1}\left( n+1\right) x\left( n+1\right) +\Delta x_{1}\left(
0\right) x\left( 0\right)   \notag \\
&&-\Delta ^{2}x_{1}\left( -1\right) x\left( 0\right) ]\frac{x_{2}\left(
n\right) }{\sqrt{\lambda \left( \lambda -4\right) }}.  \notag
\end{eqnarray}

\noindent Finally, performing necassary operations, then we get%
\begin{eqnarray}
-\sum_{i=0}^{n}\frac{q\left( i\right) x\left( i\right) x_{2}\left( i\right) 
}{\sqrt{\lambda \left( \lambda -4\right) }}x_{1}\left( n\right)
+\sum_{i=0}^{n}\frac{q\left( i\right) x\left( i\right) x_{1}\left( i\right) 
}{\sqrt{\lambda \left( \lambda -4\right) }}x_{2}\left( n\right)  &=&\frac{1}{%
\sqrt{\lambda \left( \lambda -4\right) }}[x\left( n\right) \sqrt{\lambda
\left( \lambda -4\right) }  \TCItag{3.18} \\
&&+x_{1}\left( n\right) \left( x_{2}\left( 0\right) x\left( -1\right)
-x_{2}\left( -1\right) x\left( 0\right) \right)   \notag \\
&&+x_{2}\left( n\right) \left( x_{1}\left( -1\right) x\left( 0\right)
-x_{1}\left( 0\right) x\left( -1\right) \right) ].  \notag
\end{eqnarray}%
And then, we have the following equation%
\begin{eqnarray*}
x\left( n\right)  &=&\left( \tfrac{2-h\left( 2q\left( 0\right) -2+\lambda +%
\sqrt{\lambda \left( \lambda -4\right) }\right) }{2\sqrt{\lambda \left(
\lambda -4\right) }}\right) x_{1}\left( n\right) +\left( \tfrac{-2+h\left(
2q\left( 0\right) -2+\lambda -\sqrt{\lambda \left( \lambda -4\right) }%
\right) }{2\sqrt{\lambda \left( \lambda -4\right) }}\right) x_{2}\left(
n\right)  \\
&&-\sum_{i=0}^{n}\frac{q\left( i\right) x\left( i\right) x_{2}\left(
i\right) }{\sqrt{\lambda \left( \lambda -4\right) }}x_{1}\left( n\right)
+\sum_{i=0}^{n}\frac{q\left( i\right) x\left( i\right) x_{1}\left( i\right) 
}{\sqrt{\lambda \left( \lambda -4\right) }}x_{2}\left( n\right) .
\end{eqnarray*}%
Note that, detail of the proof will be given in published paper.

\bigskip

\noindent \textbf{Theorem 3.2. }Let define Sturm-Liouville problem in
difference equations as follows;%
\begin{equation}
Ly\left( n\right) =-\lambda y\left( n\right) ,  \tag{3.19}
\end{equation}%
\begin{equation}
y\left( 0\right) =1,\text{ }y\left( 1\right) =0,  \tag{3.20}
\end{equation}%
then $\left( 3.31\right) -\left( 3.32\right) $ Sturm-Liouville problem has a
unique solution for $y\left( n\right) $ as

\begin{eqnarray}
y\left( n,\lambda \right)  &=&\left( \frac{-2+\lambda +\sqrt{\lambda \left(
\lambda -4\right) }-2q\left( 0\right) }{2\sqrt{\lambda \left( \lambda
-4\right) }}\right) \left( \frac{2-\lambda +\sqrt{\lambda \left( \lambda
-4\right) }}{2}\right) ^{n}  \TCItag{3.21} \\
&&+\left( \frac{2-\lambda +\sqrt{\lambda \left( \lambda -4\right) }+2q\left(
0\right) }{2\sqrt{\lambda \left( \lambda -4\right) }}\right) \left( \frac{%
2-\lambda -\sqrt{\lambda \left( \lambda -4\right) }}{2}\right) ^{n}  \notag
\\
&&-\sum_{i=0}^{n}\frac{q\left( i\right) y\left( i\right) \left( \frac{%
2-\lambda -\sqrt{\lambda \left( \lambda -4\right) }}{2}\right) ^{i}}{\sqrt{%
\lambda \left( \lambda -4\right) }}\left( \frac{2-\lambda +\sqrt{\lambda
\left( \lambda -4\right) }}{2}\right) ^{n}  \notag \\
&&+\sum_{i=0}^{n}\frac{q\left( i\right) y\left( i\right) \left( \frac{%
2-\lambda +\sqrt{\lambda \left( \lambda -4\right) }}{2}\right) ^{i}}{\sqrt{%
\lambda \left( \lambda -4\right) }}\left( \frac{2-\lambda -\sqrt{\lambda
\left( \lambda -4\right) }}{2}\right) ^{n}.  \notag
\end{eqnarray}

\bigskip

\noindent \textbf{Proof. }This is proved similarly to the proof of Theorem
3.1.

\bigskip 

\noindent \textbf{4.} \textbf{Asymptotic Formulas for Sturm-Liouville
Problem in Difference Equations}

\bigskip

\noindent At this section, we present the asymptotic formulas for the
solution of Sturm-Liouville problem.\textbf{\ }Let's take $\left( 3.3\right)
-\left( 3.4\right) $ Sturm-Liouville problem. Then we can give the following
Theorem.

\bigskip

\noindent \textbf{Theorem 4.1. }$\left( 3.3\right) -\left( 3.4\right) $
Sturm-Liouville problem has the estimate%
\begin{equation*}
x\left( n\right) =O\left( \mid h\mid \right) .
\end{equation*}

\bigskip 

\noindent \textbf{Theorem 4.2. }$\left( 3.19\right) -\left( 3.20\right) $
Sturm-Liouville problem has the estimate%
\begin{equation*}
y\left( n\right) =O\left( 1\right) .
\end{equation*}

\noindent The proofs will be given in published version of paper.

\bigskip 

\noindent \textbf{Conclusion}

The article has been extended the scope of the representation of solution
and asymptotic formulas for the Sturm-Liouville problem in difference
equations.

\bigskip \newpage 

\noindent \textbf{6. References}

\begin{enumerate}
\item \textbf{W. G., Kelley and A. C., Peterson,} 2001. Difference
Equations: An Introduction with Applications, Academic Press, San Diego.

\item \textbf{R. P.,} \textbf{Agarwal, }2000. Difference Equations and
Inequalities, Marcel Dekker, 970 Newyork.

\item \textbf{H., Bereketo\u{g}lu,} \textbf{V., Kutay, }2012. Fark
Denklemleri, Gazi Kitabevi, Ankara.

\item \textbf{F. V., Atkinson,} 1964. Discrete and Continuous Boundary Value
Problems, Academic Press, Newyork.

\item \textbf{A.,} \textbf{Jirari, }1995. "Second Order Sturm-Liouville
Difference Equations.and Orthogonal Polynomials, Memoirs of the American
Mathematical Society, Vol. 113, Number 542, Providence Rhode Island.

\item \textbf{B. M., Levitan, I. S., Sargsjan}, 1975. Introduction to
Spectral Theory: Selfadjoint Ordinary Differential Operators, American
Mathematical Society, Providence Rhode Island.

\item \textbf{C., M., Bender, S., A., Orszag,} 1999. Advanced Mathematical
Methods for Scirntists and Engineers: Asymptotic Methods and Perturbation
Theory, Springer-Verlag, Newyork.

\item \textbf{G., Shi, H., Wu, }2009. Spectral Theory of Sturm-Liouville
Difference Operators, Linear Algebra and Its Applications, 430, 830-846.

\item \textbf{D., B., Hinton, R., T., Lewis, }1978. Spectral Theory of
Second-Order Difference Equations, Journal of Mathematical Analysis and
Applications, 63, 421-438.

\item \textbf{Y., Shi, H., Sun, }2011.\textbf{\ }Self-adjoint extensions for
second-order symmetric linear difference equations, Linear Algebra and its
Applications, 434, 903-930.

\item \textbf{Y., Wang, Y., Shi, }2005. Eigenvalues of second-order
difference equations with periodic and antiperiodic boundary conditions,
Journal of Mathematical Analysis and Applications, 309, 56-69.

\item \textbf{S., Goldberg,} 1986. Introduction to Difference equations with
Illustrative examples from Economics, Psychology and Sociology, Dover, 260,
Newyork.

\item \textbf{V., Lakshmikantham, D., Trigiante, A. C., Peterson,} 1988.
Theory of Difference Equations: Numerical Methods and Applications,
Academic, 242, Newyork.

\item \textbf{R., Mickens,} 1990. Difference Equations, Van Nostrand, 448,
Newyork.

\item \textbf{S., Elaydi, }2005. An introduction to difference equations,
Springer Science+Business Media, Newyork.

\item \textbf{C.,} \textbf{Ahlbrandt, J., Hooker, }1987. Disconjugacy
criteria for second order linear difference equations, Qualitative
Properties of Differential Equations, Proceedings of the 1984 Edmonton
Conference, University of Alberta, Edmonton, 15-26.

\item \textbf{P., Hartman, }1978. Difference equations: disconjugacy,
principal solutions, Green's functions, Complete monocity, Trans. Amer.
Math. Soc. 246, 1-30.

\item \textbf{Y. Shi and S. Chen,} 1999. Spectral Theory of Second-Order
Vector Difference Equations, Journal of Mathematical Analysis and
Applications, 239, 195-212.

\item \textbf{S. L. Clark}, 1978. A Spectral Analysis for Self-Adjoint
Operators Generated by a Class of Second Order Difference Equations, Journal
of Mathematical Analysis and Applications, 197, 267-285.

\item \textbf{R. S. Hilscher, }2012. Spectral and oscillation theory for
general second order Sturm-Liouville difference equations, Advances in
Difference Equations, 82.

\item \textbf{J. Ji, B. Yang, }2007.\textbf{\ }Eigenvalue comparisons for
second order difference equations with Neumann boundary conditions,Linear
Algebra and its Applications, 425, 171-183.

\item \textbf{H. Sun, Y. Shi, }2006. Eigenvalues of second-order difference
equations with coupled boundary conditions, Linear Algebra and its
Applications, 414, 361-372.
\end{enumerate}

\end{document}